\documentclass[11pt]{amsart}

\usepackage{amsfonts, amssymb}
\usepackage{amscd, amsxtra, amsmath}
\usepackage{amsbsy}
\usepackage{epsfig}
\usepackage{graphicx}
\usepackage{psfrag}
\usepackage{amsmath}
\usepackage{amsthm}
\usepackage{setspace}
\usepackage{url}
\usepackage{color}
\usepackage{hyperref}
\usepackage{enumerate}
\usepackage{enumitem}

\theoremstyle{plain}
\newtheorem{theorem}{Theorem}[section]

\newtheorem{proposition}[theorem]{Proposition}

\newtheorem{question}[theorem]{Question}
\theoremstyle{definition}

\newtheorem{example}[theorem]{Example}

\newcommand{\MM}{\mathcal M}

\newcommand{\calA}{\mathcal A}

\newcommand{\BM}{\overline{\mathcal M}}

\newcommand{\OO}{\mathcal O}

\newcommand{\Eff}{\operatorname{Eff}}

\newcommand{\BEff}{\overline{\operatorname{Eff}}}

\newcommand{\bbQ}{\mathbb Q}
\newcommand{\bbR}{\mathbb R}

\newcommand{\bbP}{\mathbb P}

\newcommand{\sat}{\operatorname{sat}}

\title[Extremal cycles in moduli spaces]{An extremal effective survey about extremal effective cycles in moduli spaces of curves}

\date{\today}

\author{Dawei Chen}
\address{Department of Mathematics, Boston College, Chestnut Hill, MA 02467}
\email{dawei.chen@bc.edu}


\thanks{The author is partially supported by the NSF CAREER award DMS-1350396.}

\begin{document}

\begin{abstract}
We survey recent developments and open problems about extremal effective divisors and higher codimension cycles in moduli spaces of curves. 
\end{abstract}

\maketitle
\setcounter{tocdepth}{1}
\tableofcontents

%%%%%%%%%%%%%%%%%%%%%%%
\section{Introduction}
\label{sec:intro}
%%%%%%%%%%%%%%%%%%%%%%%

For a projective variety $X$, the cone of effective divisors $\Eff(X)$ governs its birational geometry. For instance, if the canonical divisor class of $X$ lies 
in the interior of $\Eff(X)$ (and $X$ has canonical singularities), then $X$ is of general type. Moreover, if $X$ is a Mori dream space (see \cite{HuKeel}), then $\Eff(X)$ can be decomposed into finitely many chambers such that divisor classes in the same chamber admit the same birational model of $X$ in the sense of Mori's program (or the minimal model program). In addition when $X$ is a moduli space, the resulting models often have new modular meanings.

Knowing the size of $\Eff(X)$ amounts to understanding extremal effective divisors that span the boundary rays of $\Eff(X)$. More generally, one can study the cones of effective higher codimension cycles and their extremal rays. When $X$ is a moduli space, extremal effective divisors and higher codimension cycles often parameterize objects with special geometric properties. In this survey we will focus on the case when $X$ is the 
Deligne-Mumford moduli space $\BM_{g,n}$ of stable genus $g$ curves with $n$ marked points. 

This paper is organized as follows. In Section~\ref{sec:detect} we introduce some criteria for detecting extremal effective divisors and higher codimension cycles. In Sections~\ref{sec:divisor} and~\ref{sec:higher} we review the current state of the art for extremal effective divisors and higher codimension cycles in $\BM_{g,n}$, with a focus on recent developments and open problems. To make it accessible to general audiences, we will often highlight geometric ideas and minimize technical details. Throughout the paper we work over an algebraically closed field of characteristic zero. Whenever we consider the effective cones of a variety $X$, we always assume that $X$ 
is projective and $\bbQ$-factorial. In particular, $\BM_{g,n}$ satisfies this assumption. 

\subsection*{Acknowledgements.} This survey was written at the Abel Symposium 2017. The author thanks the  organizers for the invitation and hospitality. 

%%%%%%%%%%%%%%%%%%%%%%%
\section{Detecting extremal effective cycles}
\label{sec:detect}
%%%%%%%%%%%%%%%%%%%%%%%

An effective divisor $D$ on a variety $X$ is a formal linear combination $D = \sum_{i=1}^n a_i Z_i$ where $a_i \in \bbR$ and 
$Z_i\subset X$ are subvarieties of codimension one. If all coefficients $a_i \geq 0$, we say that $D$ is effective. Denote by $\Eff(X)$ the cone of effective divisors that parameterize effective divisor classes (up to numerical equivalence $\equiv$). It is clear that $\Eff(X)$ has a convex structure. Nevertheless, 
$\Eff(X)$ may fail to be closed (see \cite[Section 1.5]{LazarsfeldPositivity} for some examples). In that case we denote by $\BEff(X)$ the closure and call it the cone of pseudoeffective divisors. 

Let $D$ be a (pseudo)effective divisor. If for any expression $D \equiv D_1 + D_2$ with $D_1$ and $D_2$ being pseudoeffective we have the divisor class of $D_i$ proportional to $D$, then $D$ is called an extremal (pseudo)effective divisor. Geometrically it means that the class of $D$ spans an extremal ray of $\BEff(X)$. Thus understanding $\BEff(X)$ amounts to finding out extremal divisors. 

How can one detect an extremal divisor? We will first illustrate the idea by a simple example, and then introduce some general criteria. 

\begin{example}
\label{ex:blowup}
Let $X$ be the blowup of $\bbP^2$ at a point $p$. Denote by $E$ the exceptional curve. Note that $E$ has two special properties. First, 
the blowdown morphism $X\to \bbP^2$ contracts $E$, i.e., $E$ is an exceptional divisor. In addition, $E^2 = -1 < 0$, i.e., $E$ has negative self intersection. 
Either one of the two properties can be used to show that $E$ spans an extremal ray of $\BEff(X)$, which we will explain below.  
\end{example}

In general for a birational contraction $f: X\dashrightarrow Y$, if $E$ is an irreducible  exceptional divisor contracted by $f$, then $E$ spans an extremal ray of $\BEff(X)$ (see e.g., \cite[Section 1.5]{RullaThesis} for a proof and generalization). In practice, it is often not easy to find a birational contraction for a given variety $X$. Alternatively, we introduce a numerical criterion which works quite well especially when $X$ is a moduli space.   

Let $D\subset X$ be an irreducible subvariety of codimension one. Suppose $C$ is an irreducible curve in $D$ such that the deformations of $C$ in $X$ cover an open dense subset of $D$. We say that $C$ is a moving curve in $D$. 

\begin{proposition}
\label{prop:divisor}
Suppose $C$ is a moving curve in $D$ such that the intersection product $C\cdot D < 0$ in $X$. Then
$D$ spans an extremal ray of $\BEff(X)$.  
\end{proposition}

To verify the above claim, suppose $D \equiv D_1 + D_2$ with $D_i$ effective and not containing $D$ in their support. Since $C\cdot D < 0$, one of $C\cdot D_i$, say $C\cdot D_1$, must be negative. It follows that $C$ (and hence any  deformation of $C$) is contained in the support of $D_1$. Since the deformations of $C$ dominate $D$, it implies that $D$ is contained in the support of $D_1$, contradicting the setup. With a little more work one can also treat the case when $D_1$ and $D_2$ are pseudoeffective (see e.g., \cite[Lemma 4.1]{ChenCoskunDivisor} for details). 

Next we turn to higher codimension cycles. Denote by $\Eff^k(X)$ (resp. $\Eff_k(X)$) the cone of effective codimension $k$ (resp. $k$-dimensional) 
cycles up to numerical equivalence, and denote by $\BEff^k(X)$ (resp. $\BEff_k(X)$) the closure. The higher codimension geometry of a variety has drawn considerable attention recently. Nevertheless, the corresponding results and techniques are still much less developed compared to the theory of divisors. 

For instance, let us try to adapt the above numerical criterion to detect extremal effective higher codimension cycles. 
Let $Z\subset X$ be a subvariety of codimension $k$. Suppose $Y\subset X$ is a $k$-dimensional subvariety 
such that $Y \cdot Z < 0$. We would like to show that $Z$ contains $Y$. However, the argument breaks down here, as the negative intersection can be caused by that $Y$ and $Z$ intersect along a positive dimensional locus, e.g., when $Y$ is a surface, $Z$ has codimension two, and $Y\cap Z$ is a curve. Therefore, we ask the following natural question. 
 
\begin{question}
Generalize Proposition~\ref{prop:divisor} to higher codimension cycles. 
\end{question} 

Now we revisit the idea of detecting extremal effective divisors from exceptional loci of birational contractions. 
Let $f: X\to Y$ be a morphism. For a 
subvariety $Z\subset X$ of dimension $k$, define the contraction index $e_f(Z) = \dim Z - \dim f(Z)$. We want to consider such $Z$ that drops maximum dimension under $f$ compared to all the other $k$-dimensional subvarieties. 

\begin{proposition}
\label{prop:higher}
Fix two integers $k > m \geq 0$. Among all $k$-dimensional subvarieties $Z$ of $X$, assume that only finitely many of them, denoted by $Z_1, \ldots, Z_n$, satisfy that $e_f(Z) \geq k - m$. If the classes of $Z_1, \ldots, Z_n$ are linearly independent, then each $Z_i$ spans an extremal ray of $\BEff_k(X)$.
\end{proposition}

For a proof and further discussion of the above result, see~\cite[Section 2]{ChenCoskunHigher}. We also state some related questions as follows. 

\begin{question}
Can one weaken the assumption of Proposition~\ref{prop:higher}, e.g., without the presence of a morphism $f$, or if $f$ is only a rational map? 
\end{question} 

\begin{question}
Does there exist a $k$-dimensional subvariety $Z\subset X$ such that $Z$ is extremal in $\Eff_k(X)$ but not extremal in $\BEff_k(X)$? 
\end{question}

%%%%%%%%%%%%%%%%%%%%%%%
\section{Extremal effective divisors in $\BM_{g,n}$}
\label{sec:divisor}
%%%%%%%%%%%%%%%%%%%%%%%

We first recall some basic facts about the moduli space $\BM_{g,n}$ of stable genus $g$ curves with $n$ (ordered) marked points. Denote by 
$\Delta = \BM_{g,n}\setminus \MM_{g,n}$ the boundary of $\BM_{g,n}$, which is a union of irreducible divisorial boundary strata $\Delta_0$ and 
$\Delta_{i; S}$ for $1\leq i \leq [g/2]$ and $S\subset \{1,\ldots, n \}$ (satisfying the stability condition). A general point of $\Delta_0$ parameterizes an irreducible curve with one node. A generic point of $\Delta_{i; S}$ parameterizes a genus $i$ curve union a genus $g-i$ curve at a node such that the marked points labeled by $S$ are contained in the genus $i$ component. These boundary divisors can further intersect, which gives a stratification of $\Delta$ based on the topological type of stable nodal curves. 

It is not hard to see that the boundary divisors are extremal in $\BEff(\BM_{g,n})$. For instance, the Torelli map $\tau$ sending a smooth curve to its Jacobian extends 
from $\BM_g$ to the Satake compactification $\calA_g^{\sat}$ of the moduli space of principally polarized Abelian varieties. The boundary divisors $\Delta_0, \ldots, \Delta_{[g/2]}$ are contracted by $\tau$, as in $\calA_g^{\sat}$ 
the Jacobian of a nodal curve does not depend on the position of the nodes. Therefore, $\Delta_0, \ldots, \Delta_{[g/2]}$ are exceptional divisors of $\tau$, hence they are extremal. Alternatively, sliding a fixed genus $g-i$ curve along a fixed genus $i$ curve gives a moving curve $C_i$ in $\Delta_i$, and  if $ i > 2$ then 
$C_i \cdot \Delta_i = 2 -2 i < 0 $, which gives another way to verify extremality using Proposition~\ref{prop:divisor}. 

It remains interesting to study whether there exist extremal effective divisors in $\BM_{g,n}$ which are not contained in the boundary. We call such divisors non-boundary extremal effective divisors. Below we will give a brief overview about their study. We separate the discussion in several cases depending on the range of $g$ and $n$. 

First, consider the case when $g$ and $n$ are both small. For instance, we know all extremal rays of $\BEff(\BM_{g,n})$ in the range $g = 0$ and $n \leq 6$; $g=1$ and $n \leq 2$; $g=2$ and $ n \leq 1$; $g=3$ and $n = 0$ (see e.g., \cite{HassettTschinkel, RullaThesis}). Even in this small range, there already exist non-boundary extremal effective divisors. For example in $\BM_3$, the locus of hyperelliptic curves forms a non-boundary extremal effective divisor. Another example is for $g = 0$. Fulton conjectured that $\BEff(\BM_{0,n})$  is generated by boundary divisors only. Keel and Vermeire  found the first counterexample on $\BM_{0,6}$ (see \cite{Vermeire}), constructed as follows. Map $\BM_{0,6}$ into the boundary of $\BM_{3}$ by gluing the six marked points in three pairs to form three nodes. Pulling back (the closure of) the locus of hyperelliptic curves thus provides a non-boundary extremal effective divisor in $\BM_{0,6}$. Unfortunately outside of this small range (even for $\BM_4$, $\BM_{2,2}$, $\BM_{1,3}$, and $\BM_{0,7}$), we do not know all extremal rays of $\BEff(\BM_{g,n})$. 

Next, consider the case when $n$ is small and $g$ is large. For example, consider the case $n= 0$, i.e., $\BM_g$ without marked points. The study of effective divisors in $\BM_g$ dates back to a series of seminal papers by Harris, Mumford, and Eisenbud (see \cite{HarrisMumfordKodaira, HarrisKodaira, EisenbudHarrisKodaira}), aiming at understanding the Kodaira dimension of $\BM_g$, which was further generalized by Farkas (see \cite{FarkasSyzygy, FarkasKoszul, FarkasSurvey}). One of the key steps is to construct effective divisors with low slopes (see \cite{HarrisMorrisonSlope}), where such divisors often parameterize curves with special linear series in the sense of Brill-Noether theory or syzygies. The Brill-Noether divisors 
are known to have the smallest slope when $g \leq 9$ and $g=11$ (see \cite{ChangRan, Tan}). For $g=10$, Farkas and Popa found an effective divisor parameterizing curves lying on $K3$ surfaces whose slope is smaller than the Brill-Noether divisors (see \cite{FarkasPopa}). We refer to \cite{ChenFarkasMorrison} for a detailed introduction and further discussion along this circle of ideas. Despite these results, for general $g$ we do not know any non-boundary extremal effective divisor in $\BM_g$. 

\begin{question}
Find non-boundary extremal effective divisors in $\BM_g$ for general $g$. 
\end{question}

One can also consider curves with marked points, say $\BM_{g,1}$ with one marking only. One of the most natural effective divisors in $\BM_{g,1}$ is the Weierstrass divisor $W$ parameterizing curves with a marked Weierstrass point. For small $g$ we know that $W$ is an extremal effective divisor (see \cite{RullaThesis, Jensen56, Jensen34, ChenCycle}). 

\begin{question}
Is $W$ extremal in $\BEff(\BM_{g,1})$ for general $g$? 
\end{question}

Note that if $W$ is an exceptional divisor of a birational contraction of $\BM_{g,1}$, then the extremality would follow right away (see \cite{Polishchuk} for an attempt). 

Now consider the case when $g$ is small and $n$ is arbitrary. Start with $\BM_{0,n}$ first. Using some combinatorial  graphs, Castravet and Tevelev constructed a series of non-boundary extremal effective divisors in $\BM_{0,n}$ for $n \geq 7$, called hypertree divisors, which significantly generalize the construction of Keel and Vermeire for $\BM_{0,6}$ (see \cite{CastravetTevelevHypertree}). Another related breakthrough they showed is that $\BM_{0,n}$ is not a Mori dream space for $n \geq 134$ (see \cite{CastravetTevelevMori}). The bound was then reduced to $n\geq 13$ by Gonz\'{a}lez and Karu (see \cite{GonzalezKaru}) and further to $n\geq 10$ by Hausen, Keicher, and Laface (see \cite{HKL}). Note that for any given $n$, the hypertree divisors in $\BM_{0,n}$ are finitely many. In general, being not a Mori dream space does not imply that the effective cone is not finitely generated. Therefore for any given $\BM_{g,n}$, one can still ask whether there exist only finitely many non-boundary extremal effective divisors, i.e., whether $\BEff(\BM_{g,n})$ is finitely generated. 

Joint with Coskun we showed that the above question has a negative answer for $\BM_{1,n}$, by constructing infinitely many non-boundary extremal effective divisors in $\BM_{1,n}$ for each $n\geq 3$ as follows (see \cite{ChenCoskunDivisor}). Fix an $n$-tuple $\underline{a} = (a_1, \ldots, a_n)$ of integers such that 
$\sum_{i=1}^n a_i = 0$. Define a divisor $D_{\underline{a}}\subset \BM_{1,n}$ as the closure of the locus of $(E, p_1, \ldots, p_n)$ where 
$E$ is a smooth genus one curve and $\sum_{i=1}^n a_i p_i \sim 0$. If in addition $\gcd (a_1, \ldots, a_n) = 1$, then $D_{\underline{a}}$ is irreducible. Under these conditions we showed that $D_{\underline{a}}$ spans an extremal ray of $\BEff(\BM_{1,n})$ by finding a negative moving curve in $D_{\underline{a}}$. Varying the values of $a_i$  thus provides infinitely many such extremal effective divisors. One may further ask whether these $D_{\underline{a}}$ are the only non-boundary extremal effective divisors in $\BM_{1,n}$. Indeed it is not the case, as joint with Patel we showed that there exist other non-boundary extremal effective divisors in $\BM_{1,n}$ by pulling back divisors of Brill-Noether and Gieseker-Petri type from $\BM_{g,n}$ 
(see \cite{ChenPatel}). Similarly by pulling back some of the $D_{\underline{a}}$ divisors from $\BM_{1,n-2}$ to $\BM_{0,n}$, Opie showed that there exist other non-boundary extremal effective divisors in $\BM_{0,n}$ besides the hypertree divisors (see \cite{Opie}). 
New extremal effective divisors in $\BM_{0,n}$ for $n\geq 9$ were further discovered by Doran, Giansiracusa, and Jensen under a more general study of the Cox ring of $\BM_{0,n}$ (see \cite{DGJ}). 

Finally when $g$ and $n$ are both large, we understand very little about the structure of $\BEff(\BM_{g,n})$. Among the few known results, Farkas and Verra showed that the Brill-Noether divisor in $\BM_{g,g}$ parameterizing 
$(C, p_1, \ldots, p_g)$ such that $h^0(C, p_1 + \cdots + p_g) \geq 2$ is extremal, since it is contracted under 
the birational map from $\BM_{g,g}$ to the universal Picard variety over $\BM_{g}$ given by $(C, p_1, \ldots, p_g) \mapsto \OO_C (p_1 + \cdots + p_g)$ (see \cite{FarkasVerra}). In addition, the divisors $D_{\underline{a}}$ on $\BM_{1,n}$ can also be interpreted as parameterizing $(E, p_1, \ldots, p_n)$ where $\sum_{i=1}^n a_i p_i$ is a (meromorphic) canonical divisor (i.e., a trivial divisor on a genus one curve). Using the strata of canonical divisors with a given type of zeros and poles for general $g$, recently Mullane constructed infinitely many non-boundary extremal effective divisors in each $\BM_{g,n}$ in the range $g\geq 2$ and $n \geq g+1$ (see \cite{Mullane}). 

We summarize some related questions as follows. 

\begin{question}
Find new non-boundary extremal effective divisors in $\BM_{g,n}$. In particular, is $\BEff(\BM_{0,n})$ finitely generated when $n \geq 7$ and is $\BEff(\BM_{g,n})$ finitely generated when $n \leq g$? 
\end{question}

%%%%%%%%%%%%%%%%%%%%%%%
\section{Extremal higher codimension cycles in $\BM_{g,n}$}
\label{sec:higher}
%%%%%%%%%%%%%%%%%%%%%%%

We have seen that the boundary divisors are extremal in $\BM_{g,n}$. Similarly one can ask whether the higher codimension boundary strata are extremal. Joint with Coskun we showed that 
every codimension two boundary stratum of $\BM_g$ and of $\BM_{0,n}$ is extremal, and every codimension $k$ boundary stratum of $\BM_{0,n}$ parameterizing curves with $k$ marked tails attached to an unmarked $\bbP^1$ is extremal (see \cite{ChenCoskunHigher}). The proofs use the criterion in Proposition~\ref{prop:higher} and we illustrate the idea by the following example.

Recall the extended Torelli map $\tau: \BM_g \to \calA_g^{\sat}$. Take an irreducible codimension two boundary stratum $\Delta_{ij}$ whose general point parameterizes a genus $g-i-j$ curve attached to two tails of genus $i$ and $j$, respectively. Moreover suppose that $i$, $j$, and $g-i-j $ are at least two. Then the contraction index $e_{\tau}(\Delta_{ij}) = 4$, as the position of the two nodes of a general curve parametrized by $\Delta_{ij}$ is forgotten under $\tau$. On the other hand, if $Z\subset \BM_g$ is a subvariety of codimension two, then the condition 
$e_\tau (Z) \geq 4$ implies that a general curve parameterized by $Z$ contains at least two nodes, hence $Z$ must be one of the codimension two boundary strata. Note that there are only finitely many codimension two boundary strata and their classes are linearly independent (see \cite{Edidin}). Therefore, Proposition~\ref{prop:higher} implies that $\Delta_{ij}$ is extremal in $\BEff^2(\BM_{g,n})$. 

One can also look for non-boundary extremal higher codimension cycles in $\BM_{g,n}$. In \cite{ChenCoskunHigher} we showed that the locus of hyperelliptic curves with a marked Weierstrass point is a non-boundary extremal codimension two cycle in $\BM_{3,1}$ and the locus of hyperelliptic curves is a non-boundary extremal codimension two cycle in $\BM_4$. The proofs of these results are more complicated compared to the cases of divisors and boundary strata, which reply on some delicate analysis of canonical curves in genus three and four. Moreover, joint with Tarasca we showed that the locus of genus two curves with $n$ marked Weierstrass points is an extremal codimension $n$ cycle in $\BM_{2,n}$ for $n\leq 6$ (see \cite{ChenTarasca}). Blankers further showed that the locus of hyperelliptic curves with $\ell$ marked Weierstrass points, $m$ marked conjugate pairs of points and $n$ free marked points is an extremal codimension $(\ell+m)$ cycle in 
$\BM_{2, \ell+2m+n}$ (see \cite{Blankers}). 

Even for small values of $g$ and $n$, the effective cone of higher codimension cycles in $\BM_{g,n}$ can be very complicated. For instance, Schaffler constructed more extremal codimension two cycles in $\BM_{0,7}$ besides the boundary strata and the lifts of the Keel-Vermeire divisors from $\BM_{0,6}$ (see \cite{Schaffler}). In \cite{ChenCoskunHigher} we showed that there exist infinitely many extremal effective codimension two cycles in $\BM_{1,n}$ for every $n \geq 5$ and in $\BM_{2,n}$ for every $n \geq 2$, respectively. Using the strata of canonical divisors, recently Mullane generalized the idea of \cite{ChenCoskunHigher} and exhibited infinitely many extremal effective codimension $k$ cycles in $\BM_{g,n}$ for every $g$ and $n$ in the range $g\geq 3$, $n \geq g-1$ and $k=2$; $g\geq 2$, $k\leq n-g, g$; $g=1$, $k\leq n-2$ (see \cite{MullaneHigher}).  

We end the survey by the following questions. 

\begin{question}
Is $\BEff^2(\BM_{0,7})$ finitely generated? 
\end{question}

\begin{question}
Is the locus of hyperelliptic curves in $\BM_g$ extremal for $g\geq 5$? 
\end{question}

%%%%%%%%%%%%%%%%%%%%%%%

\end{document}